\documentclass{amsart}
\usepackage{amsmath, amssymb, amsthm}
\usepackage{mathtools}
\usepackage{amscd}
\usepackage{graphicx} 
\usepackage{todonotes}
\usepackage{tabularx}
\usepackage{nicefrac}
\usepackage{framed}
\usepackage{xcolor}
\colorlet{shadecolor}{green!25}
\usepackage{hyperref}
\usepackage{url}
\usepackage{enumitem}

\newtheorem{thm}{Theorem}[section]

\newtheorem{prop}[thm]{Proposition}
\newtheorem{lemma}[thm]{Lemma}

\newtheorem{cor}[thm]{Corollary}

\newtheorem{Polyunits}[thm]{Auxiliary polynomials for units}
\newtheorem{Polyint}[thm]{Auxiliary polynomials for integers}

\newtheorem*{thm*}{Theorem}
\newtheorem*{alg*}{Algorithm}
\newtheorem*{lemma*}{Lemma}

\theoremstyle{remark}

\newtheorem*{rmk*}{Remark}

\newtheorem*{notation*}{Notation}

\newtheorem*{example*}{Example}

\theoremstyle{definition}

\newtheorem*{defn*}{Definition}

\newcommand{\mybf}{\mathbb}

\newcommand{\Q}{\mybf{Q}}
\newcommand{\Z}{\mybf{Z}}
\newcommand{\C}{\mybf{C}}

\author[Pottmeyer]{Lukas Pottmeyer}
\email{lukas.pottmeyer@uni-due.de}
\address{Universit\"{a}t Duisberg-Essen, 45117 Essen, Germany}

\title{Weighted error terms in lower bounds for totally $p$-adic integers}

\date{\today}

\begin{document}

\begin{abstract}
We introduce a method for obtaining lower bounds for the heights of totally $p$-adic algebraic integers and units using weighted combinations of auxiliary polynomials. By combining the associated archimedean error terms, we obtain better estimates than those arising from using individual auxiliary polynomials. Our main application of these weighted height bounds is the calculation of a sharp lower bound for the height of totally $2$-adic algebraic integers, namely $\frac{1}{2}\log 2$. Moreover, we show that this minimum is isolated: every such integer, apart from the elements of height zero and the four minimizers, has height at least $\frac{8}{15}\log 2$.
\end{abstract}

\maketitle 
\author


\section{Introduction}

In \cite{Sch} Schinzel proved that the (absolute logarithmic Weil-)height $h$ of a totally real algebraic integer different from $0$ and $\pm1$ is bounded from below by $\frac{1}{2}\log\left(\frac{1+\sqrt{5}}{2}\right)$. Moreover, this bound is sharp since it is attained for $\frac{1+\sqrt{5}}{2}$. By a classical result of Kronecker \cite{Kr}, the elements of height zero are precisely $0$ and the roots of unity. Hence, we have to exclude these elements in the search for lower height bounds.

As a first version of a $p$-adic analogue of Schinzel's result, Bombieri and Zannier calculated in \cite{BZ} an effective lower bound for the $\liminf$ of the heights of elements $\alpha$, whose Galois orbit can be embedded in a fixed finite extension of $\Q_p$ for prime numbers $p$. We call an $\alpha\in\overline{\mathbb{Q}}$ totally $p$-adic, if the Galois orbit of $\alpha$ can be embedded in $\Q_p$. Hence, we may call a totally real element \emph{totally $\infty$-adic}. Denote for every place $v$ on $\Q$ the field of totally $v$-adic numbers with $\mathbb{Q}^{tv}$. Let $M$ be a set of rational primes, and define $F_M=\cap_{p\in M} \Q^{tp}$. Restricting to totally $v$-adic numbers, the result of Bombieri and Zannier reads
\[
\liminf_{\alpha \in F_M} h(\alpha) \geq \frac{1}{2}\sum_{p\in M} \frac{\log(p)}{p+1}.
\] 
This result was refined in \cite{Fi},\cite{FiPe}, and \cite{FiPo} and allows also to include $\infty\in M$. Lower bounds for the height of totally $p$-adic algebraic units have been calculated by Petsche \cite{Pe}, and Dubickas and Mossinghoff \cite{DM}. In \cite{Po} the author extended these results to lower bounds for totally $p$-adic algebraic units, integers and numbers, respectively. Define
\begin{align*}
u_p &= \inf\{h(\alpha) \vert \alpha \text{ an algebraic unit } \in \Q^{tp}\setminus\{\text{roots of }1\}\} \\
i_p &= \inf\{h(\alpha) \vert \alpha \text{ an algebraic integer } \in (\Q^{tp})^*\setminus\{\text{roots of }1\}\} \\
n_p &= \inf\{h(\alpha) \vert \alpha  \in (\Q^{tp})^*\setminus\{\text{roots of }1\}\} .
\end{align*} 

The main result of \cite{Po} was

\begin{thm}\label{thm:LP}
The following bounds hold

\bigskip

\begin{minipage}{0.4\textwidth}
For an odd prime $p$ we have
\begin{align*}
\frac{\log(p/2)}{p-1} < &u_p < \frac{\log(p+\nicefrac{1}{p^{p-2}})}{p-1}\\
\frac{\log(p/2)}{p} < &i_p \leq \frac{\log(p)}{p} \\
\frac{\log(p/2)}{p+1} < &n_p .
\end{align*}
\end{minipage} \qquad
\begin{minipage}{0.4\textwidth}
In the case $p=2$ we have
\begin{align*}
\log(2) < &u_2 < \frac{1}{2}\log\left(\frac{4+\sqrt{17}}{2}\right)\\
\frac{2}{5}\log(2)< &i_2 \leq \frac{1}{2}\log(2) \\
\frac{1}{4}\log(2) < &n_p .
\end{align*}
\end{minipage}
\end{thm}

The error term $\log(2)$ in these lower bounds comes solely from the archimedean valuation, which leads to the estimate $h(\alpha-1)\leq h(\alpha)+\log(2)$. In the case of function fields (where no archimedean valuations exists), the same reasoning as in the proof of  Theorem \ref{thm:LP} leads to sharp lower bounds, as demonstrated in \cite{FaPe}. In this note, we will present a method for shrinking this archimedean error term by using a weighted error term for several polynomials instead of only $f(x)=x-1$. The idea is quite simple: Choose integer polynomials $f_1(x)$ and $f_2(x)$, such that we can control the $p$-adic valuation of $f_1(\alpha)$ and $f_2(\alpha)$ for all totally $p$-adic algebraic integers, resp. units. Moreover, these polynomials should satisfy that the error $h(f_1(\alpha))-\deg(f_1)h(\alpha)$ is small, whenever the error $h(f_2(\alpha))-\deg(f_2)h(\alpha)$ is big. To make this precise we define for any $f\in\mathbb{Z}[x]$ and $z\in\mathbb{C}$ the value
\[
c_f(z) = \max\{0,\log\vert f(z)\vert \}-\deg(f)\max\{0,\log(\vert z\vert\}.
\]
Let $f_1,\ldots,f_k \in\mathbb{Z}[x]\setminus\mathbb{Z}$, and take a set of real numbers $t_1,\ldots,t_k\geq 0$ such that $\sum_{i=1}^k t_k =1$. Set $\boldsymbol f=(f_1,\ldots,f_k)$ and $\boldsymbol t=(t_1,\ldots,t_k)$. Let $S$ be a subset of $\mathbb{Q}_p$ that does not contain any root of the polynomials $f_1,\ldots,f_k$. Let $v_p$ a valuation on $\Q_p$ normalized such that $v_p(p)=1$. We define the weighted $p$-adic contribution as
\begin{equation}\label{eq:defm}
m_p(S,\boldsymbol{f},\boldsymbol{t})
=
\inf_{x\in S}
\sum_{i=1}^k
\frac{t_i}{\deg(f_i)}
\max\{0,v_p(f_i(x))\},
\end{equation}
and the weighted error term as
\begin{equation}\label{eq:defE}
E(\boldsymbol{f},\boldsymbol{t}) = \sup_{z\in\C,\vert z\vert =1} \sum_{i=1}^k \frac{t_i}{\deg(f_i)}c_{f_i}(z).
\end{equation}

\begin{lemma}\label{mainlemma}
We use the notation from above. Let $\alpha$ be an algebraic number, such that all Galois conjugates of $\alpha$ lie in $S$. Then we have
\[
h(\alpha)\geq m_p(S,\boldsymbol{f},\boldsymbol{t})\log(p) - E(\boldsymbol{f},\boldsymbol{t}).
\]
\end{lemma}

This is the main lemma of this note. As a first application of this lemma, we will choose for an odd prime $p$ the two polynomials $f_1(x)=x^{p-1}-1$ (as usual) and $f_2(x)=\frac{p-1}{2}x^{p-1}+\frac{p+1}{2}$. Note that the maximum of $c_{f_1}(z)$ is attained at for $z^{p-1}=-1$ and that $z^{p-1}=-1$ gives $c_{f_2}(z)=0$. We can calculate by hand the best choices for the weights. This gives 

\begin{prop}\label{prop:up}
For any odd prime $p$ we have $u_p \geq \frac{\log\left(\frac{\sqrt{p^2+3}}{2}\right)}{p-1}$. 
\end{prop}

Lemma \ref{mainlemma} should be compared to the height estimates of Dubickas and Mossinghoff \cite{DM}. There only units and only the contribution of  one polynomial were considered. Moreover, they used the value 
\[
\frac{\log(\omega_p(f))}{\deg(f)} = \frac{1}{\deg(f)}\log\left(\gcd\left(
f(1),\, pf'(1),\, \frac{p^2f''(1)}{2},\,
\ldots,\, \frac{p^{\deg(f)}f^{(\deg(f))}(1)}{\deg(f)!}
\right)\right)
\]
 instead of our $m_p(S,\boldsymbol{f},\boldsymbol{t})\log(p)$. The simple example of $k=1$, $f_1(x)=x^2 -1$ and $p=2$ shows, that $m_p(S,\boldsymbol{f},\boldsymbol{t})\log(p)$ may be strictly larger than $\frac{\log(\omega_p(f))}{\deg(f)}$.  Dubickas and Mossinghoff found an iterative algorithm which produces a sequence of polynomials for which their lower bound for $u_p$ increases. The degree of these polynomials rapidly increases as well. The first  iteration of the algorithm in \cite{DM} produces the general lower bound $u_p \geq \frac{\log\left(\frac{\sqrt{p^2+1}}{2}\right)}{p-1}$.

In contrast to this iterative algorithm, we will present a finite family of polynomials with degree $\leq 2$ in $x^{p-1}$. If we choose some of these polynomials, the best choice of the parameters $t_1,\ldots,t_k$ apparently can be found by linear programming, as already applied by Smyth in \cite{Sm}. \footnote{ChatGPT (GPT-6) wrote a Python code which seems to produce correct results. The code is available on the authors webpage \cite{LPurl}. The results presented in this note are all proven rigorously within this note, independently of the numerical optimization of the weights.} 

For small primes, our set of auxiliary functions is of moderate size. Searching through subsets of these polynomials leads to much better bounds for $u_p$ (and $i_p$) in the cases $p\in\{2,3,5\}$. Dubickas and Mossinghoff raised the question for which primes we have $n_p<u_p$. This inequality is already known for $p=2$ \cite{DM} and $p=3$ \cite{Po}. Our method gives also:

\begin{prop}\label{prop:u5}
It is $u_5 > 0{,}26923 >\frac{\log(5)}{6}\geq n_5$.
\end{prop}

Numerical values for lower bounds for $u_2$ and $u_3$ will be presented in the last section. For $p=7$ we could not find an application of Lemma \ref{mainlemma} proving $n_p<u_p$.

We now turn to lower bounds for $i_p$. In the case $p=2$ we see that the elementary bounds $\frac{2}{5}\log(2) < i_2 \leq \frac{1}{2}\log(2)$ from Theorem \ref{thm:LP} are quite narrow. We found auxiliary polynomials in our prescribed list, for which Lemma \ref{mainlemma} closes the gap and allows us to calculate the value of $i_p$ explicitly.

\begin{thm}\label{thm:i2}
We have $i_{2}=\frac{1}{2}\log(2)$. Moreover, we have $h(\alpha) \geq \frac{8}{15}\log(2)$ for all totally $2$-adic algebraic integers $\alpha\notin \{0,\pm 1, \pm \frac{1\pm\sqrt{7}i}{2}\}$.
\end{thm}

In \cite{St} and \cite{AD} sharp lower bounds for the height of totally $p$-adic algebraic numbers of degree $3$ and $4$ have been calculated. Our bound is sharp among totally $2$-adic algebraic integers of any degree. To the best of my knowledge, this is the first known sharp lower bound for totally $p$-adic integers. It remains open whether the equation $n_2=i_2$ holds.

The paper is organized as follows. In Section \ref{sec:mainlemma} a proof of Lemma \ref{mainlemma} is given, and the lists of auxiliary polynomials are presented. Section \ref{sec:general} contains -- as an (hopefully) illustrative example -- the proof of Proposition \ref{prop:up}. Proposition \ref{prop:u5} is proven in Section \ref{sec:u5}. The sharp lower bound for the heights of totally $2$-adic algebraic integers from Theorem \ref{thm:i2} is proven in Section \ref{sec:maintheorem}. Finally, in Section \ref{sec:numerical} we present several numerical result for lower bounds for $u_p$ and $i_p$ for $p\in \{2,3,5\}$.

\section{The main lemma}\label{sec:mainlemma}

Let $f\in\Z[x]\setminus\Z$ be a polynomial, and $\alpha$ be a totally $p$-adic algebraic integer, with $f(\alpha)\neq 0$.  Let $S$ be any subset of $\Q_p$ not containing any root of $f$, such that all conjugates of $\alpha$ lie in $S$. Then we have
\begin{align}\label{eq:heightfpos}
h(f(\alpha))=h((f(\alpha))^{-1})&\geq \frac{\log(p)}{[\mathbb{Q}(\alpha):\Q]} \sum_{v\mid p} \inf_{x\in S} \max\{0,v_p(f(x))\} \nonumber \\ &= \log(p)  \inf_{x\in S} \max\{0,v_p(f(x))\}.
\end{align}
Hence, the quantity $\inf_{x\in S} \max\{0,v_p(f_i(x))\}$ controls the height of $h(f(\alpha))$. We need to compare this with $h(\alpha)$. Recall that for any $f\in \mathbb{Z}$ we defined
\[
c_f(z) = \max\{0,\log\vert f(z)\vert \}-\deg(f)\max\{0,\log\vert z\vert\}.
\]
We have
\begin{equation}\label{eq:heightdifference}
h(\alpha) \geq \frac{1}{\deg(f)} \cdot\left(h(f(\alpha))- \frac{1}{[\mathbb{Q}(\alpha):\Q]} \sum_{\sigma} c_f(\sigma(\alpha))\right),
\end{equation}
where $\sigma$ runs over all embeddings $\Q(\alpha)\longrightarrow \overline{\Q}$. See Lemma 3.1 in \cite{Po}.

As in the introduction let $f_1,\ldots,f_k \in\mathbb{Z}[x]\setminus\mathbb{Z}$, and take a set of real numbers $t_1,\ldots,t_k\geq 0$ such that $\sum_{i=1}^k t_k =1$. Set $\boldsymbol f=(f_1,\ldots,f_k)$ and $\boldsymbol t=(t_1,\ldots,t_k)$. Lastly, let $S$ be any subset of $\mathbb{Q}_p$ that does not contain any root of the polynomials $f_1,\ldots,f_k$. We defined
\begin{equation*}
m_p(S,\boldsymbol{f},\boldsymbol{t})
=
\inf_{x\in S}
\sum_{i=1}^k
\frac{t_i}{\deg(f_i)}
\max\{0,v_p(f_i(x))\},
\end{equation*}

and the weighted error term
\[
E(\boldsymbol{f},\boldsymbol{t}) = \sup_{z\in\C,\vert z\vert =1} \sum_{i=1}^k \frac{t_i}{\deg(f_i)}c_{f_i}(z).
\]

Hence, when one $c_{f_i}(z)$ is large, this may be compensated by a small value $c_{f_j}(z)$ for some of the other polynomials. By \eqref{eq:heightdifference} and \eqref{eq:heightfpos} we have

\begin{align*}
h(\alpha) &=\sum_{i=1}^k t_i h(\alpha) \geq \sum_{i=1}^k t_i \frac{1}{\deg(f_i)} \cdot\left(h(f_i(\alpha))- \frac{1}{[\mathbb{Q}(\alpha):\Q]} \sum_{\sigma} c_{f_i}(\sigma(\alpha))\right)\\
&\geq \sum_{i=1}^k t_i \frac{1}{\deg(f_i)} \cdot\left(h(f_i(\alpha))\right) - \sup_{z\in\C}  \sum_{i=1}^k \frac{t_i}{\deg(f_i)}c_{f_i}(z)\\
&\geq m_p(S,\boldsymbol{f},\boldsymbol{t})\log(p) - E(\boldsymbol{f},\boldsymbol{t}).
\end{align*}

The last inequality follows from the definition of $m_p(S,\boldsymbol{f},\boldsymbol{t})$, and from the maximum modulus principle for subharmonic functions. This proves Lemma \ref{mainlemma}.

In the case $S=\mathbb{Q}_p \setminus \{z\in \Q_p\vert f_1(z)\cdots f_k(z)=0\}$, this will give a trivial bound since in this case $m_p(S,\boldsymbol{f},\boldsymbol{t})=0$. This can be seen by calculating $v_p(f_i(p^{-n}))$ for increasing $n$. This obstacle may be removed  by using $h(\alpha)=h(\alpha^{-1})$. However, the lower bounds we achieve this way, are too weak to justify this additional technicality. Hence, we are mainly interested in the choices $S_p(\boldsymbol{f})=\mathbb{Z}_p \setminus \{z\in \Q_p\vert f_1(z)\cdots f_k(z)=0\}$, which calculates lower bounds for totally $p$-adic integers, and $S_p^*(\boldsymbol{f})=\mathbb{Z}_p^* \setminus \{z\in \Q_p\vert f_1(z)\cdots f_k(z)=0\}$, which calculates lower bounds for totally $p$-adic algebraic units. Here, we denote with $\Z_p$ the ring of  $p$-adic integers.

In the case of elements $\alpha$ which are simultaneously totally $p$-adic and totally real, the bound can be strengthened further. The only difference is the weighted error-term, which can be replaced by

\[
E_{\mathbb{R}}(\boldsymbol{f},\boldsymbol{t}) = \sup_{z\in \mathbb{R}} \sum_{i=1}^k \frac{t_i}{\deg(f_i)}c_{f_i}(z).
\]

\begin{lemma}\label{mainlemmareal}
We use the notation from above. Let $\alpha$ be an algebraic number, such that all Galois conjugates of $\alpha$ lie in $S$. Then we have
\[
h(\alpha)\geq m_p(S,\boldsymbol{f},\boldsymbol{t})\log(p) - E_{\mathbb{R}}(\boldsymbol{f},\boldsymbol{t}).
\]
\end{lemma}

 We close this section by presenting sets of potentially auxiliary polynomials. For any element  $x\in \Z_p^*$, we use $v_p(x^{p-1}-1)\geq 1$. Hence, we only consider polynomials in the variable $x^{p-1}$.
 
 \begin{Polyunits}\label{polyunits}
The set of auxiliary polynomials in the search for lower height bounds for totally $p$-adic algebraic units consists of:
\begin{itemize}
\item $x^{p-1} -1$
\item $ax^{p-1} + b$, with integers $a,b \in \{1,\ldots, p\}$ satisfying $a+b=p$
\item $a(x^{p-1})^2 + bx^{p-1} + c$, with integers $a,b,c$, satisfying $p^2 \mid a+b+c$, $p\mid 2a+b$, and $\vert a \vert ,\vert b\vert, \vert c\vert \leq p^2$
\item any products of different polynomials from above 
\end{itemize}
These polynomials $f$ are chosen such that $v_p(f(x))$  is at least the degree of $f$ as a polynomial in $y=x^{p-1}$ for every $p$-adic unit $x$.
\end{Polyunits}

f we are looking for auxiliary polynomials to prove lower bounds for the height of totally $p$-adic algebraic integers, we need to include more polynomials in order to control the $p$-adic valuation for elements $x\in \Q_p$ satisfying $v_p(x)>0$.

\begin{Polyint}\label{polyint}
The set of auxiliary polynomials in the search for lower height bounds for totally $p$-adic algebraic integers consists of:
\begin{itemize}
\item $x$
\item $ax^2+bx+c$, with $\vert a \vert ,\vert b\vert, \vert c\vert \leq p^2$.
\item $x^{p-1} -1$
\item  $ax^{p-1} + b$, with integers $a,b \in \{1,\ldots, p\}$ satisfying $a+b=p$
\item $a(x^{p-1})^2 + bx^{p-1} + c$, with integers $a,b,c$, satisfying $p^2 \mid a+b+c$, $p\mid 2a+b$, and $\vert a \vert ,\vert b\vert, \vert c\vert \leq p^2$
\item $a(x^{p-1})^2 + bx^{p-1} + c$, with integers $a,b,c$, satisfying $p^2 \mid a+b+c$ and $\vert a \vert ,\vert b\vert, \vert c\vert \leq p^2$
\item any products of different polynomials from above 
\end{itemize}
\end{Polyint}

\section{General lower bounds for totally $p$-adic units}\label{sec:general}

Let $p$ be an odd prime. We define $f_1(x)=x^{p-1}-1$, $f_2(x)=\frac{p-1}{2}x^{p-1}+\frac{p+1}{2}$, and $\boldsymbol{f}=(f_1,f_2)$. For any $t\in[0,1]$ we define $\boldsymbol{t}=(t_1,t_2)=(t,1-t)$. Now, let $\alpha$ be a totally $p$-adic algebraic unit, which is not a root of unity. Then $f_1(\alpha)\neq 0\neq f_2(\alpha)$. We apply Lemma \ref{mainlemma}, with $S$ equal to $p$-adic units different from roots of unity. Note, that whenever $c_{f_1}(z)$ is close to $2$, we have $c_{f_2}(z)$ close to zero. 

By Fermat's little theorem we have $m_p(S,\boldsymbol{f},\boldsymbol{t})=\frac{1}{p-1}$, since $\deg(f_1)=\deg(f_2)=p-1$. It remains to bound the weighted error term from below. We have
\begin{align*}
E(\boldsymbol{f},\boldsymbol{t})&=\sup_{z\in \C,\vert z\vert =1} \left( \frac{t}{p-1} (\max\{0,\log\vert f_1(z)\vert \}-(p-1)\max\{0,\log\vert z\vert\})\right. \\ & \qquad \left.+ \frac{1-t}{p-1}( \max\{0,\log\vert f_2(z)\vert \}-(p-1)\max\{0,\log\vert z\vert\})\right) \\
&=\frac{1}{p-1} \cdot \sup_{z\in \C,\vert z\vert =1}  \left( t\max\{0,\log\vert z^{p-1}-1\vert\} \right.\\ & \qquad \left.+(1-t)\max\{0,\log\vert  \frac{p-1}{2}z^{p-1}+\frac{p+1}{2}\vert \} \right).
\end{align*}

Since on the unit circle we have $\vert z^{p-1}\vert =\vert z\vert =1$, we may write $z^{p-1}=e^{i\theta}$ and calculate
\[
\sup_{\theta\in[0,2\pi]}  \left( t\max\{0,\log\vert e^{i\theta}-1\vert\}+(1-t)\max\{0,\log\vert  \frac{p-1}{2}e^{i\theta}+\frac{p+1}{2}\vert \} \right).
\]
Writing $e^{i\theta}=\cos(\theta)+i\sin(\theta)$, we have $\vert e^{i\theta}-1\vert =\sqrt{2-2s}$ and $\vert  \frac{p-1}{2}e^{i\theta}+\frac{p+1}{2}\vert= \sqrt{\frac{p^2-1}{2} s + \frac{p^2+1}{2}}$, with the substitution $s=\cos(\theta)$.  This leads to

\[
(p-1) E(\boldsymbol{f},\boldsymbol{t})= \sup_{s\in[-1,1]} \begin{cases} \frac{t}{2}\log(2-2s) +  \frac{1-t}{2}\log\left(\frac{p^2-1}{2}s + \frac{p^2+1}{2}\right)  &\text{ if } -1 \leq s < \frac{1}{2} \\ 
									\frac{1-t}{2}\log\left(\frac{p^2-1}{2} s+ \frac{p^2+1}{2}\right)  &\text{ if } \frac{1}{2} \leq s \leq 1 \end{cases}
\]

Denote with $s_t$ a critical value of the right hand side, where the maximum is attained. We want to pick a $t$ such that the maximum of the right hand side is minimal. Hence,  $s_t \in [-1,\frac{1}{2}]$, because otherwise the choice of $t=1$ would give $E(\boldsymbol{f},\boldsymbol{t})=0$. We have $s_t=1- t\cdot \frac{2p^2}{p^2-1}$, which lies in the given interval whenever $t\in[\frac{p^2-1}{4p^2},\frac{p^2-1}{p^2}]$. Hence, for every $t\in[\frac{p^2-1}{4p^2},\frac{p^2-1}{p^2}]$ we have
\[
(p-1) E(\boldsymbol{f},\boldsymbol{t})=\frac{t}{2}\log\left(t\frac{4p^2}{p^2-1}\right) + \frac{1-t}{2}\log(p^2(1-t)).
\]
The only critical value of this function is $t=\frac{p^2-1}{p^2+3}$, which lies indeed in the forced interval for admissible values of $t$. Now we can fix $\boldsymbol{t}=(  \frac{p^2-1}{p^2+3} , 1- \frac{p^2-1}{p^2+3})$, and with this particular choice we finally get
\[
(p-1) E(\boldsymbol{f},\boldsymbol{t}) = \frac{1}{2} \log\left( 4\cdot \frac{p^2}{p^2+3} \right) = \log\left( 2 \sqrt{\frac{p^2}{p^2+3}} \right) .
\]
Hence, the error term was indeed pushed below $\log(2)$. Lemma \ref{mainlemma} now gives
\[
h(\alpha)\geq \frac{\log(p) - \log\left( 2 \sqrt{\frac{p^2}{p^2+3}} \right)}{p-1} =\frac{\log\left( \frac{1}{2} \sqrt{p^2+3} \right)}{p-1},
\]
which proves $u_p\geq \frac{\log\left(\frac{\sqrt{p^2+3}}{2}\right)}{p-1}$ and hence Proposition \ref{prop:u5}.

\begin{cor}\label{cor:npip}
For every odd prime $p$, we have
\[
i_p\geq \frac{\log\left(\frac{\sqrt{p^2+3}}{2}\right)}{p} \quad \text{ and } n_p\geq \frac{\log\left(\frac{\sqrt{p^2+3}}{2}\right)}{p+1}.
\]
\end{cor}
\begin{proof}
For an arbitrary element $\alpha\in \mathbb{Q}^{tp}$ we may define  $r$ as the number of extensions $v\mid p$ in $\Q(\alpha)$, with $v(\alpha)>0$, and similarly $s$ as the number of extensions $v\mid p$, with $v(\alpha)<0$. With $h(\alpha)=h(\alpha^{-1})$ we may assume $r\geq s$. Moreover $h(\alpha)\geq \frac{r}{[\Q(\alpha):\Q]}\log(p)$, and the results from above imply
\[
h(\alpha)\geq \frac{[\Q(\alpha):\Q]-r-s}{[\Q(\alpha):\Q]} \frac{\log(p)}{p-1}- \frac{1}{p-1}\log\left( 2 \sqrt{\frac{p^2}{p^2+3}} \right).
\]
By choosing $r$ and $s\leq r$ such that the maximum of theses two bounds in minimized, we find the claimed bounds for $n_p$ and $i_p$ (where $s=0$).
\end{proof}

\section{Proof of Proposition \ref{prop:u5}}\label{sec:u5}

For any choice of polynomials from \ref{polyunits} in $\boldsymbol{f}=(f_1,\ldots,f_k)$ we can calculate the optimal vector $\boldsymbol{t}$ in Lemma \ref{mainlemma} by linear programming. Our search leads to  $k=2$, 
\[
f_1(x)=x^4 -1 \quad \text{ and } \quad f_2(x)=(x^4-1)(2x^4+3)(3x^4+2)(9x^8 +7x^4 +9) .
\]
The numerically optimal weights are claimed to be $t'_1 =  0{,}4213\ldots$ and $t'_2=1-t_1$. In order to get rigorous results, we approximate this and use $\boldsymbol{t}=(\frac{3}{7},\frac{4}{7})$.

\begin{lemma}\label{A1}
With the notation as above, we have 
\[
m_5(S_5^*(\boldsymbol{f}),\boldsymbol{f},\boldsymbol{t})\geq \frac{39}{140}.
\]
\end{lemma}
\begin{proof}
By Fermat, we know  that for every $x\in S_p^*(\boldsymbol{f})$ we have $x^4=1+p\cdot u$, for some $u\in \Z_p$. .
It follows immediately that  $v_5(f_1(x))\geq 1$ and $v_5(f_2(x))\geq 5$ for every $x\in S_p^*(\boldsymbol{f})$. Distinguishing between the possible residue classes of $u$ modulo $p$, we see that at least one factor of $f_2(x)$ gives an additional $+1$ to the $p$-adic valuations. Hence, $v_5(f_2(x))\geq 6$. It follows 
\[
m_p(S_p^*(\boldsymbol{f}),\boldsymbol{f},\boldsymbol{t})
\geq \frac{3}{7}\cdot\frac{1}{4} + \frac{4}{7}\cdot \frac{6}{20}=\frac{39}{140}.
\]   
\end{proof}

\begin{lemma}\label{A2}
We have $E(\boldsymbol{f},\boldsymbol{t})=\frac{1}{280} \log\left( (2-2\tau)^{19} (12\tau+13)^8(18\tau+7)^8 \right)$, where $\tau =-\frac{513}{1512}-\sqrt{\frac{18875}{84672}}$. This is:  $E(\boldsymbol{f},\boldsymbol{t})= 0{,}1791072\ldots$
\end{lemma}
\begin{proof}
For $\vert z \vert =1$ we have $\vert z^4\vert =1$. Write $z^4=e^{i\tau}$ with real part $s=\cos(\tau)$. Then 
$\vert f_1(z)\vert^2 =(2-2s)$ and $\vert f_2(z)\vert^2=(2-2s)(12s+13)^2(18s+7)^2$. Hence, we need to maximize the function $g:[-1,1]\longrightarrow \mathbb{R} $, with
\[
g(s) = \frac{1}{2} \left(\frac{3}{7\cdot 4}\max\{0,\log(2-2s)\} + \frac{4}{7\cdot 20} \max\{0,\log((2-2s)(12s+13)^2(18s+7)^2)\}\right).
\]
We separate the cases $s\in [-1, \frac{1}{2}]$ and $s\in(\frac{1}{2},1]$. In both cases the extremal values can be calculated by hand, and we find that the maximum of $g$ is attained at $s=-\frac{513}{1512}-\sqrt{\frac{18875}{84672}}$. 
Plugging this back in gives the claimed value.
\end{proof}

Now, Proposition \ref{prop:u5} follows easily. Let $\alpha$ be a totally $5$-adic algebraic unit, not a root of unity. Then by Lemmas \ref{A1}, \ref{A2}, and \ref{mainlemma} we have
\[
h(\alpha)\geq \frac{39}{140}\log(5) - 0{,}17911 =0{,}269233\ldots
\]
The statement $n_5\leq \frac{\log(5)}{6}$  is proven in Corollary 5.4 of \cite{Po}.

\section{Totally $2$-adic algebraic integers}\label{sec:maintheorem}

Since we consider algebraic integers instead of units, we take polynomials from \ref{polyint}. We found the three polynomials
\[
f_1(x)=x\quad , \quad f_2(x)=x^2-1 \quad , \quad  f_3(x)=x^4 + 3x^2 +4.
\]
Note that we need $f_1(x)$ to control the height for elements with positive $2$-adic valuation. The polynomial $f_2(x)=(x-1)(x+1)$ was already used in the proof of \ref{thm:LP}, and $f_3(x)=(x^2 + x+2)(x^2-x+2)$ vanishes precisely at the smallest totally $2$-adic algebraic integers we could find. Hence, these three candidates are quite natural candidates from \ref{polyint}. The python code gives the approximated weights $t_1=\frac{11}{30}$, $t_2=\frac{3}{10}$, and $t_3=\frac{1}{3}$. Hence, we apply Lemma \ref{mainlemma}, with $\boldsymbol{f}=(f_1,f_2,f_3)$, and $\boldsymbol{t}=(t_1,t_2,t_3)$.

We start by calculating 
\begin{align*}
 & m_2(S_2(\boldsymbol{f}),\boldsymbol{f},\boldsymbol{t}) \\ = &\inf_{x\in S_2(\boldsymbol{f})}\left(\frac{11}{30}  \max\{0,v_2(x)\} + \frac{3}{20} \max\{0,v_2(x^2 -1)\} + \frac{1}{12} \max\{0,v_2(x^4+3x^2+4)\} \right).
\end{align*}

If $v_2(x)=0$, then -- since the square of any any odd number is congruent to $1$ modulo $8$ -- we write $x^2=1+2^3u$, with $u\in \Z_2$. We see that $v_p(f_2(x))\geq 3$ and $v_p(f_3(x))\geq 3$. Distinguishing between $v_2(u)=0$ and $v_2(u)>0$, at least one of the terms $f_2(x),f_3(x)$ gives an additional $+1$ to the valuation $v_2$. Of course $v_2(f_1(x))=0$ in this case. 

If on the other hand $v_2(x)>0$, we write $x=2u$, with $u\in\Z_p$. It follows $v_2(f_3(x))=  v_2(2^4 u^4 + 3\cdot 2^2 u^2 + 2^2)=v_2(2^2(2^2 + 3u^2+1))$. Hence, if $v_2(u)>0$ we have $v_2(f_3(x))=2$. If however $v_2(u)=0$ and we write $u=1+2u'$, we get $v_2(f_3(x))=(2^4 (3u'(u'+1)+2))\geq 5$. Of course $v_2(f_2(x))=0$.

Hence, we have four cases depending on $v_2(x)=0$ and $v_2(u)=0$ respectively. Taking the least favorable of these cases gives
\begin{equation}\label{eq:mi2}
m_2(S_2(\boldsymbol{f}),\boldsymbol{f},\boldsymbol{t})\geq \frac{47}{60}.
\end{equation}
 
Before we calculate the weighted error term, we note that by definition $c_{f_1}=0$. Hence, we 
 \begin{align*}
 E(\boldsymbol{f},\boldsymbol{t}) &= \sup_{z\in\C,\vert z\vert =1}\left( \frac{3}{20} \max\{0,\log\vert z^2-1\vert\} + \frac{1}{12} \max\{0,\log\vert z^4+3z^2+4\vert\}\right)  \\
 \end{align*}
 For any $z\in\C$, with $\vert z\vert =1$, set $z^2=e^{i\theta}$ and $s=\cos(\theta)$. Then $\vert f_1(z)\vert^2 = 2-2s$ and $\vert f_2(z)\vert = 16s^2 + 30 s +18$. Hence
 \[
 E(\boldsymbol{f},\boldsymbol{t}) =\sup_{s\in[-1,1]} \left( \frac{3}{40} \max\{0,\log(2-2s)\} + \frac{1}{24} \max\{0,\log\vert 16s^2+30s+18\vert\}\right).
 \]
 Again, we separately calculate the maximum of this function for $s\in[-1,\frac{1}{2}]$ and for $s\in(\frac{1}{2},1]$. In the first case, the derivative of the right hand side is $\frac{76s^2+65s+3}{480s^3+420s^2-360s-540}$. This leads to the critical values $s=-1$, $s=\frac{-65\pm\sqrt{3313}}{152}$, and $s=\frac{1}{2}$. In the second case we note that only the last summand appears, which takes its maximum at $s=1$. Comparing these possible values, gives that the supremum is attained at $s=1$. Hence, we achieve the nice value
 \begin{equation}\label{eq:Ei2}
  E(\boldsymbol{f},\boldsymbol{t})  = \frac{1}{24} \log(64)=\frac{1}{4}\log(2).
 \end{equation}
 Note that this is considerably smaller than the error term $\log(2)$ from the elementary bound.
 
 Now we can finally prove Theorem \ref{thm:i2}. For any totally $2$-adic integer not a root of $f_1(x)\cdot f_2(x)\cdot f_3(x)$, Lemma \ref{mainlemma} together with \eqref{eq:mi2} and \eqref{eq:Ei2}, yields
 \[
 h(\alpha)\geq \log(2)m_2(S_2(\boldsymbol{f}),\boldsymbol{f},\boldsymbol{t}) - E(\boldsymbol{f},\boldsymbol{t}) \geq \log(2) \frac{47}{60}.-\frac{1}{4}\log(2)=\frac{8}{15}\log(2).
 \]
 The roots of $f_1(x)\cdot f_2(x)\cdot f_3(x)$ we had to exclude, are $0$, $\pm1$, and the four numbers $\pm \frac{1\pm\sqrt{7}i}{2}$. Each of these for last numbers is totally $2$-adic and has height $\frac{1}{2}\log(2)<\frac{8}{15}\log(2)$. Hence, $i_2=\frac{1}{2}\log(2)$. Moreover, the value $\frac{1}{2}\log(2)$ is isolated in the spectrum of heights of totally $2$-adic algebraic integers. This is the same phenomenon that appears in the field of totally real numbers, where the minimal positive value $\frac{1}{2}\log(\frac{1+\sqrt{5}}{2})$ is isolated, as shown by Smyth \cite{Sm}.
 
Also in \cite{BZ} it is shown that $\liminf_{\alpha \in \Z_2}h(\alpha)\leq \log(2)$. (See also \cite{CeFe}, \cite{Fi}, and \cite{PeSt} for generalizations and different proofs of this result.) Hence, it would be interesting to determine the number of totally $p$-adic algebraic integers of height between $\frac{8}{15}\log(2)$ and $\log(2)$.

\section{Further numerical values}\label{sec:numerical}

We will list lower bounds for $u_2$, $u_3$, and $i_3$, $i_5$. These were produced by the following lists of auxiliary polynomials, with the given weights. The results were computed numerically. A check with exact arithmetic calculations has not been done.

\begin{itemize}
\item[$u_2$:] In the calculation of $u_2$ we increased the possible size of the coefficients in the list \ref{polyint} and also allowed polynomials of degree up to $4$. The choice 
\begin{align*}
f_1(x)&=x^2-1  \\f_2(x)&=17x^4+30x^2+17 \\ f_3(x)&=23x^4+18x^2+23,
\end{align*}
and $\boldsymbol t=\frac{1}{1000}\left(582,237,181\right)$, gives $u_2 >0{,}7402$.
\item[$u_3$:] The choice
\begin{align*}
f_1&=7x^4+4x^2+7\\
f_2&=5x^4-x^2+5\\
f_3&=(x^2-1)(2x^2+1)\\
f_4&=(x^2-1)(x^2+2),
\end{align*}
and $\boldsymbol t=\frac18(1,1,3,3)$, gives $u_3>0{,}4467$.
\item[$i_3$:] The choice 
\begin{align*}
f_1&=x(x^2-1)\\
f_2&=x(x^2+2)\\
f_3&=x(2x^2+1)\\
f_4&=x,
\end{align*}
and $\boldsymbol t=\frac1{12}(5,3,3,1)$, gives $i_3>0{,}2441$.
\item[$i_5$:] The choice
\begin{align*}
f_1&=x(x^4-1)\\
f_2&=x(9x^8+7x^4+9)\\
f_3&=x(2x^4+3)\\
f_4&=x(3x^4+2)\\
f_5&=x\\
f_6&=x(x^4-1)(3x^4+2),
\end{align*}
and $\boldsymbol t=
\frac1{10000}(5072,2036,1131,954,488,319)$, gives $i_5>0,2107$.
\end{itemize}

It remains open, if this method can potentially prove other sharp lower bounds as $i_2$. To this end, the sets of auxiliary polynomials need to be increased, which makes the search for good candidates infeasible. An algorithmic approach which generates for any prime $p$ good choices for the auxiliary functions $f_1,\ldots,f_k$ is still missing.

\end{document}